# A Universal Simulation Platform for Flexible Systems

**Minh Hoang Vu,**[1,*] **Tajwar Abrar Aleef,**[1] **Usama Pervaiz,**[1] **Yeman Brhane Hagos**[1] **and Saed Khawaldeh**[1]

[1] Erasmus+ Joint Master Program in Medical Imaging and Applications University of Burgundy (France), University of Cassino (Italy) and University of Girona (Spain)

**Abstract.** This article proposes a universal simulation platform for simulating systems undergoing duress. In other words, this paper introduces a total simulation package which includes a number of methods of simulating the flexibility of a given system. This platform includes detailed procedures for simulating a flexible link by a numerical method called finite difference method. In order to verify the effectiveness of the proposed process, two examples are covered in different situations to discuss the importance of boundary control and mesh selection in the way of ensuring the stability of the system. In addition, a graphical user interface (GUI) application called the SimuFlex is designed having a selection of methods that the user can choose along with the parameters of the controllers that can be easily manipulated from the GUI.

## 1 Introduction

Flexible systems including strings and beams are widely applied in many areas of modern mechanical engineering [1], such as aerospace, civil [2], marine [3], chemical engineering [4], and especially in automatic manufacturing assembly. Systems with lightweight flexible links are becoming more popular as because they possess many advantages over the conventional rigid-link ones such as better speed, lower cost, greater labor productivity, enhanced machine efficiency, improved quality, increased system reliability, and reduced parts inventories to name a few. With the aim of improving industrial productivity, it is desirable to build systems with lightweight materials to increase the payload-to-weight ratio. This article is motivated by the industrial applications in the simulation of vibration of flexible structures.

**Corresponding author:**
Minh Hoang Vu, E-mail: hoangminh.vu@smartdatics.com

Simulation is incredibly important as it helps to identify and solve problems before implementation. Furthermore, it can be used to predict the courses and results of certain actions, thus it enables to test hypotheses without having to carry them out in the practical world. Choosing a simulation method plays a key role in the successful simulation. Finite difference method (FDM) was chosen due to its simplicity and straightforward implementation. There are not many fast solvers and packages that are using FDM algorithm for partial difference equations (PDEs). Additionally, although Matlab provides a powerful visualization workspace, it does not have the functionality to allow researchers to input specifications and controller design at the same time. Engineers and researchers normally need to use a combination of tools for model building, controller design, and numerical simulation. Solving these difficulties greatly motivates the research in designing a GUI application, written on Matlab, for a complete system capable of simulating flexible-link systems.

For the purpose of dynamic modeling and simulation, flexible systems are regarded as distributed parameter systems which are in infinite dimension with various boundary conditions involving functions of space and time. In practice, dynamics of flexible risers are mathematically represented by a set of PDEs with appropriate boundary equations or approximated by ordinary differential equations (ODEs). Since the vibration of flexible links are governed by PDEs, flexible systems are distributed-parameter and possess an infinite number of dimensions which makes them difficult to control[5]. Vibration reduction to minimize the disturbance effects is desirable for preventing damage and improving lifespan in industry. To avoid the vibration problem in a system, some popular methods, such as the application of point-care maps [6], variable structure control [7], sliding model control [8], energy-based robust control , model-free control [9], [10], distributed model predictive control [11], and boundary control [12] [13] are considered. In this article, both examples use boundary control.

In this article, by using a well-understood numerical method, namely FDM, we design the boundary control based on the distributed parameter model of the flexible system. Then the importance of choosing appropriate control parameters can be discovered to improve the outcome of the final result. The main contribution of this article includes presenting a procedure for some kind of flexible systems represented by PDEs. Next, two examples are provided in order to demonstrate the procedure above. Finally, a GUI application, SimuFlex, for flexible systems is introduced.

The rest of this article is organized as follows. Section 2 develops a procedure for simulating a flexible system using FDM. Section 3 shows two applications of the proposed process. Next, in Section 4, a GUI application, written by MATLAB for simulation of some specific flexible systems, is developed. Finally, Section 5 summarises and concludes the article.

## 2 Procedure for simulating a flexible system using finite difference method

### 2.1 Step 1: Modeling problem using Hamilton's method

With Lagrangian mechanics, Hamilton's equations provide a new and equivalent way of looking at classical mechanics. These equations, however, do not provide a convenient way of solving particular problems. In addition, these equations provide deeper insights into both the general structure of classical mechanics and its connection to quantum mechanics as understood through Hamiltonian mechanics as well as its connection to remaining areas of science.

Using Hamilton's equations by following these steps [14]:

(i) First write out the Lagrangian $L = T - V$. Express $T$ and $V$ as though you were going to use Lagrange's equation.

(ii) Calculate the momenta by differentiating the Lagrangian with respect to velocity:

$$p_i(q_i, \dot{q}_i, t) = \frac{p\mathcal{L}}{p\dot{q}_i} \quad (1)$$

(iii) Express the velocities in terms of the momenta by inverting the expressions in step (ii).

(iv) Compute the Hamiltonian using the usual definition of H as the *Legendre transformation* of L:

$$\mathcal{H} = \sum_i q_i \frac{p\mathcal{L}}{p\dot{q}_i} - \mathcal{L} = \sum_i \dot{q}_i p_i - \mathcal{L} \quad (2)$$

and then substitute for velocity in step (iii)

(v) Apply Hamilton's equations.

Most likely, in step 1, researcher is expected to bring out governing equation, boundary initial conditions, disturbance, and controller of the system.

### 2.2 Step 2: Select Mesh Configuration

Generating a grid that is a finite set of points, which the beam is divided into. Normally, element configuration can be selected by choosing a specific number of mesh in time or space. In this article, $N$ and $T$ are space and time grids, respectively. Space and time grids have to be chosen appropriately but not too small or large because these selections might affect seriously the stability of the system. There must be a relationship between some crucial parameters, which assure the stability. This problem is discussed in the last part of this section.

### 2.3 Step 3: Apply boundary and initial conditions

Applying boundary and initial conditions in general equation in order to find a solution uniquely, which is invoked from Keller's theorem [15]. These conditions are commonly specified in terms of known values of the unknowns on a part of the surface.

### 2.4 Step 4: Define general equation for each element

Substituting the derivatives with some finite difference formulas into the general equation for each point, thus we can obtain the algebraic system of equations. This method can be called the procedure of imposing the Dirichlet BCs [16]. Some important finite difference approximation equations are listed as follows [17]:

$$w_t(i,j) = \frac{w(i, j+1) - w(i,j)}{k} + O(k)$$

$$w_{tt}(i,j) = \frac{w(i, j+1) - 2w(i,j) + w(i, j-1)}{k^2} + O(k^2)$$

$$w_i(i,j) = \frac{w(i+1, j) - w(i,j)}{h} + O(h)$$

$$w_{xx}(i,j) = \frac{w(i+1, j) - 2w(i,j) + w(i-1, j)}{h^2} + O(h^2)$$

$$w_{xxx}(i,j) = \frac{w(i+2, j) - 2w(i+1, j) + 2w(i-1, j) - w(i-2, j)}{h^3} + O(h^3)$$

where $h = L/N$, and $k = t_f/T$, while $L$ is the length of the flexible link and $t_f$ is the time we run simulation on. Furthermore, $w_t(i,j)$, and $w_x(i,j)$ represent the first derivatives of $w$ with respect to $t$ and $x$, respectively, at $x = i$ and $t = j$. Hence, $w_{tt}(i,j)$ is the second derivative with respect to $t$, $w_{ttt}(i,j)$ is the third derivative with respect to $t$, and so on.

### 2.5 Step 5: Solve for derived or secondary quantities by recursion

Programming the problem to obtain the approximate solution. In specific, in each time node, a recursion is built with the results of boundary and initial conditions from *step 3*, and governing equation from *step 4* in order to compute the displacement of each node along the flexible link.

### 2.6 Step 6: Interpretation the result

Analyzing the stability of the result. When the result can not converge, we need to go back to *step 1* to check carefully the controller following by the appropriateness of $T$ and $N$ in *step 2* and start again.

Consequently, error analysis plays the key considerations in the development and application of all numerical methods. While it is an extensively cultivated field, the tools available are often inadequate, especially in nonlinear problems [18].

How to choose mesh placements $N$ and time step $T$ to ensure the consistency, convergence, and stability?

#### 2.6.1 Scheme for the Heat Equation

$$\phi_t(i,j) = \alpha \phi_{xx}(i,j) \tag{3}$$

Finite difference approximation gives

$$\frac{\phi(i,j+1) - \phi(i,j)}{k} = \frac{\phi(i-1,j) - 2\phi(i,j) + \phi(i+1,j)}{h^2} + \mathcal{O}(k) + \mathcal{O}(h^2) \tag{4}$$

that is

$$\phi(i,j+1) = \phi(i,j) + \alpha k/h^2 \Big[\phi(i-1,j) - 2\phi(i,j) + \phi(i+1,j)\Big] \tag{5}$$

or

$$\phi(i,j+1) = r\phi(i+1,j) + (1-2r)\phi(i,j) + r\phi(i-1,j) \tag{6}$$

Stable solution with the Heat scheme is obtained only if and only if:

$$r = \alpha k/h^2 < 1/2 \tag{7}$$

#### 2.6.2 Scheme for the Flexible Beam Equation

Governing equation of the system is:

$$\rho w_{tt}(i,j) = -EI w_{xxxx}(i,j) + T w_{xx}(i,j) - c w_t(i,j) + f(i,j) \tag{8}$$

for all $(i,j) \in (0,N) \times [0,\infty)$.

Approximating each components using FD and grouping them together gives:

$$\begin{aligned}
w(i,j) =& 2w(i,j-1) - w(i,j-2) \\
&+ k^2T/\rho h^2 \Big[w(i+1,j-1) - 2w(i,j-1) \\
&+ w(i-1,j-1)\Big] - k^2EI/\rho h^4 \Big[w(i+2,j-1) \\
&- 4w(i+1,j-1) + 6w(i,j-1) \\
&- 4w(i-1,j-1) + w(i,j)\Big] \\
&- kc/\rho \Big[w(i,j-1) - w(i,j-2)\Big] \\
&+ k^2/\rho f(i,j-1)
\end{aligned} \tag{9}$$

Stable condition holds when:

$$4k^2 EI/\rho h^4 + k^2 T/\rho h^2 \leq 1 \tag{10}$$

## 3 Two problems

Two examples are investigated to demonstrate the above procedure for simulating a flexible system in this section. To discover the importance of appropriate controllers and mesh configurations that are $N$ and $T$, plan was to choose appropriate controller and mesh configurations in example 1 and 2, respectively.

### 3.1 Timoshenko Beam

#### 3.1.1 Step 1: Modeling problem using Hamilton's method

Consider a Timoshenko beam in Figure 1. Denote the displacement of the Timoshenko beam at position $i$, and time $j$ by $w(i,j)$, and rotation of the Timoshenko beam's cross-section owing by $\varphi(i,j)$.

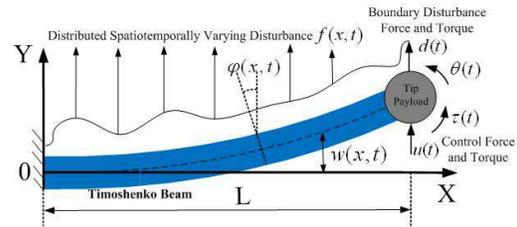

**Figure 1.** Timoshenko beam

By modeling this flexible system, we obtain governing, boundary equations, and initial conditions as follows. Governing equations of T-beam is defined as:

$$\rho w_{tt}(i,j) = -K\Big[\varphi_x(i,j) + w_{xx}(i,j)\Big] + f(i,j) \quad (11)$$

$$I_\rho \varphi_{tt}(i,j) = EI\varphi_{xx}(i,j) - K\Big[\varphi(i,j) + w_x(i,j)\Big] \quad (12)$$

for all $(i,j) \in (0,N) \times [0,\infty)$, where $\rho$ is uniform mass per unit length of the Timoshenko beam.

$K = kGA$ is a positive constant that depends on the shape of the Timoshenko beam's cross-section. $I_\rho$ is the uniform mass moment while inertia of the Timoshenko beam's cross-section and $EI$ is bending stiffness of the Timoshenko beam.

The disturbance $d(i,j)$ and $\theta(i,j)$ on the tip payload is generated by:

$$\begin{aligned} d(i,j) &= 1 + \sin(\pi ij) + \sin(2\pi ij) + \sin(3\pi ij) \\ \theta(i,j) &= \sin(\pi ij) + \sin(2\pi ij) + \sin(3\pi ij) \end{aligned} \quad (13)$$

Boundary conditions of the system are written as:

$$w(0,j) = \varphi(0,j) = 0 \quad (14)$$

$$Mw_{tt}(N,j) - K\Big[\varphi(N,j) - w_x(N,j)\Big] = u(j) + d(i,j) \quad (15)$$

$$J\varphi_{tt}(N,j) + EI\varphi_x(N,j) = \tau(j) + \theta(i,j) \quad (16)$$

where $w(N,j)$ is the tip point of the beam. $M$ denotes mass of the payload. $J$ is the inertia of the payload. $u(j)$ and $\tau(j)$ are the controllers.

Next, the corresponding initial conditions of the beam system are given as

$$w(i,0) = i/2 \quad (17)$$

$$\varphi(i,0) = \pi/6 \quad (18)$$

$$\varphi_t(i,0) = w_t(i,0) = 0 \quad (19)$$

Finally, distributed disturbance along the beam is described

$$\begin{aligned} f(i,j) = i/{1000L}\Big[&1 + \sin(0.1\pi ij) \\ &+ \sin(0.2\pi ij) + \sin(0.3\pi ij)\Big] \end{aligned} \quad (20)$$

#### 3.1.2 Step 2: Select Mesh Configuration

Choose $T = 10000$ and $N = 50$.

#### 3.1.3 Step 3: Apply boundary and initial conditions

Initial conditions approximation from Eq. (19) and noting Eq. (17) and Eq. (18)

$$w(i,1) = w(i,0) = i/2 \quad (21)$$
$$\varphi(i,1) = \varphi(i,0) = \pi/6 \quad (22)$$

Boundary condition approximation from Eq. (14) means that the first node is fixed

$$w(0,j) = \varphi(0,j) = 0 \quad (23)$$

In this article, two cases are investigated, comprising without control and with PD control.

**(i) Without control**

$$u(j) = 0 \quad (24)$$

$$\tau(j) = 0 \quad (25)$$

Substituting Eq. (24) into Eq. (15) to obtain

$$\begin{aligned} w(N,j) =& 2w(N,j-1) - w(N,j-2) \\ &+ {k^2K}/{M}\Big[\varphi(N,j-1)\Big] + {k^2}/{M}\Big[d(N,j)\Big] \\ &- {k^2K}/{Mh}\Big[w(N,j-1) \\ &- w(N-1,j-1)\Big] \end{aligned} \quad (26)$$

Substituting Eq. (25) into Eq. (16) to obtain

$$\begin{aligned} \varphi(N,j) =& 2\varphi(N,j-1) - \varphi(N,j-2) \\ &- {k^2EI}/{Jh}\Big[\varphi(N,j-1) - \varphi(N-1,j-1)\Big] \\ &+ {k^2}/{J}\Big[\theta(N,j)\Big] \end{aligned} \quad (27)$$

**Remark.** *The graph of Timoshenko beam without control is shown in Figure 2. As we observe, the tip position of the beam fluctuates around 0, which reflects the bad effect of external sources or undesired disturbance. However, the system still keeps its stability.*

**(ii) PD control**

$$u(j) = -k_1 w(N,j) - k_1 w_t(N,j) \quad (28)$$

$$\tau(j) = -k_3 \varphi(N,j) - k_4 \varphi_t(N,j) \quad (29)$$

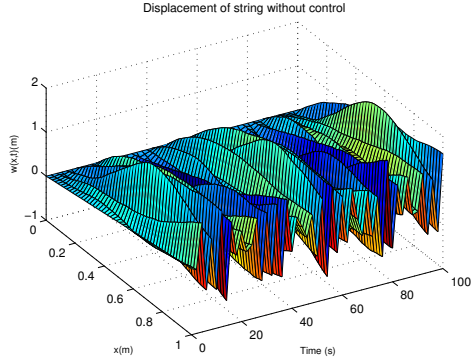

**Figure 2.** T beam - without control

Substituting Eq. (28) into Eq. (15)

$$\begin{aligned}w(N,j) =& 2w(N,j-1) - w(N,j-2) \\&+ k^2K/M\Big[\varphi(N,j-1)\Big] + k^2/M\Big[d(N,j)\Big] \\&- k^2K/Mh\Big[w(N,j-1)\Big] \\&- k^2k_d/M\Big[w(N,j-1)\Big] \\&- k^2k_p/Mk\Big[w(N,j-1) - w(N,j-2)\Big] \\&- w(N-1,j-1)\Big]\end{aligned}$$
(30)

Substituting Eq. (29) into Eq. (16)

$$\begin{aligned}\varphi(N,j) =& 2\varphi(N,j-1) - \varphi(N,j-2) \\&- k^2EI/Jh\Big[\varphi(N,j-1) - \varphi(N-1,j-1)\Big] \\&- k^2k_d/J\Big[\varphi(N,j-1)\Big] \\&- kk_p/J\Big[\varphi(N,j-1) - \varphi(N,j-2)\Big] \\&+ k^2/J\Big[\theta(N,j)\Big]\end{aligned}$$
(31)

#### 3.1.4 Step 4: Define general equation for each element

Finally, we obtain the approximation for governing equation from Eq. (11)

$$\begin{aligned}w(i,j) =& w(i,j-1) - w(i,j-2) + k^2/\rho\Big[f(i,j)\Big] \\&- k^2K/\rho h\Big[\varphi(i,j-1) - \varphi(i-1,j-1)\Big] \\&+ k^2/\rho h^2\Big[w(i+1,j-1) \\&- 2w(i,j-1) + w(i-1,j-1)\Big]\end{aligned}$$
(32)

And Eq. (12)

$$\begin{aligned}\varphi(i,j) =& 2\varphi(i,j-1) - \varphi(i,j-2) \\&+ k^2EI/I_ph^2\Big[\varphi(i+1,j-1) - 2\varphi(i,j-1) \\&+ \varphi(i-1,j-1)\Big] - Kk^2/I_p\Big[\varphi(i,j-1)\Big] \\&+ k^2/I_ph\Big[w(i,j-1) - w(i-1,j-1)\Big]\end{aligned}$$
(33)

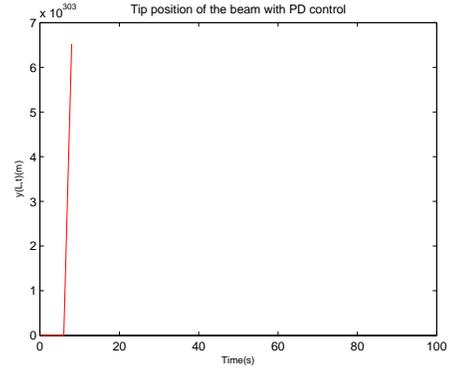

**Figure 3.** Tbeam - PD control ($k_1 = k_3 = 100$ and $k_2 = k_4 = 30$)

**Remark.** *Unsuitable controllers $k_1 = k_3 = 100$ and $k_2 = k_4 = 30$, which make the system lose stability, is shown in Figure 3. The shape displacement goes up dramatically without decreasing a trend. Thus, we can conclude that large controllers lead the system to be unstable.*

#### 3.1.5 Step 6: Interpretation the result

As we can observe that large controllers $k_1 = k_3 = 100$ and $k_2 = k_4 = 30$ makes the beam lose its stability, to avoid this problem, we might reduce the magnitude of $k_2$ to 10.

**Remark.** *Figure 4 and Figure 5 reflect the beam with PD control while we choose suitable controllers $k_1 = k_3 = 100$ and $k_2 = k_4 = 10$. The shape displacement sharply decreases from initial condition 1 to about 0. Then, it flutters around this point for a short period before converging as expected.*

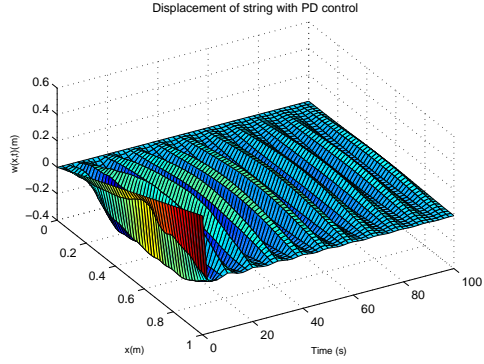

**Figure 4.** Tbeam - PD control ($k_1 = k_3 = 100$ and $k_2 = k_4 = 10$)

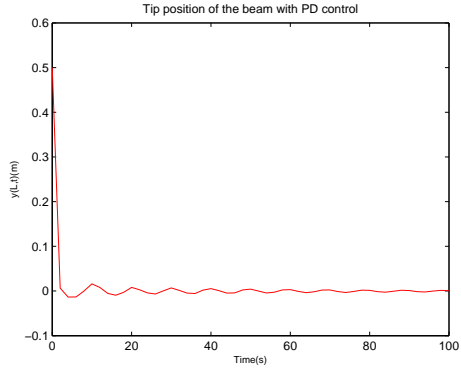

**Figure 5.** T beam - PD control ($k_1 = k_3 = 100$ and $k_2 = k_4 = 10$)

### 3.2 Non-uniform String

#### 3.2.1 Step 1: Modeling problem using Hamilton's method

Consider a non-uniform string in Figure 6. Denote the displacement of the non-uniform string at position $i$, and time $j$ by $w(i, j)$.

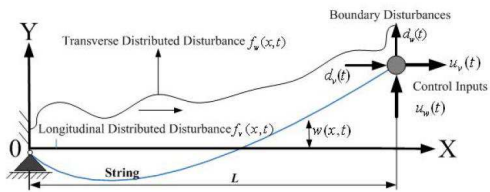

**Figure 6.** Non-uniform string

Governing equation of the non-uniform string system is represented by

$$\rho(i)w_{tt}(i,j) = T(i,j)w_{xx}(i,j) + T_x(i,j)w_x(i,j) \\ + \lambda_x(i)[w_x(i,j)]^3 \\ + 3\lambda(i)[w_x(i,j)]^2 w_{xx}(i,j) + f(i,j) \quad (34)$$

for all $(i, j) \in (0, N) \times [0, \infty)$, where $T(i, j)$ denotes the tension at position $x = i$, and time $t = j$. $\rho(i)$ is uniform mass per unit length, which is dependent on the position of the mesh.

It can be observed that $T$ and $\rho$ in this case are different from the ones in the above mentioned example, where they are constants. The tension $T(i, j)$ and $\lambda(i)$ of the string can be expressed as

$$T(i,j) = T_0(i) + \lambda(i)[w_x(i,j)]^2 \quad (35)$$

where $T_0(i) = 10(i+1)$ and $\lambda(i) = 0.1i$. Boundary conditions of the system is given

$$w(0, j) = 0 \quad (36)$$

$$Mw_{tt}(N,j) + T(N,j)w_x(N,j) + \lambda(L)[w_x(N,j)]^3 \\ = u(t) + d(t) \quad (37)$$

The corresponding initial conditions of the string system are given as

$$w(i, 0) = i \quad (38)$$

$$w_t(i, 0) = 0 \quad (39)$$

For simulation study, the boundary disturbance $d(j)$ on the tip payload is generated by the following equation

$$d(j) = 1 + 0.2\sin(0.2j) + 0.3\sin(0.3j) + 0.5\sin(0.5j) \quad (40)$$

The time-varying distributed disturbance $f(i, j)$ on the string is described as

$$f(i,j) = i \times \Big[3 + \sin(\pi i j) + \sin(2\pi i j) + \sin(3\pi i j)\Big] \quad (41)$$

#### 3.2.2 Step 2: Select Mesh Configuration

Choose $T = 10000$ and $N = 50$.

#### 3.2.3 Step 3: Apply boundary and initial conditions

Initial conditions approximation from Eq. (39)

$$w(i, 0) = w(i, 1) \quad (42)$$

Boundary condition approximation for Eq. (36)

$$w(0, j) = 0 \tag{43}$$

In this article, two cases are simulated, comprising without control and with with exact-model control.

**(i) Without control**

$$u(j) = 0 \tag{44}$$

The approximation for Eq. (37) when $u(j) = 0$

$$\begin{aligned}w(N, j) =\ & 2w(N, j-1) - w(N, j-2) \\ & - {}^{k^2 T(N,j)}\!/\!{}_{Mh}\Big[w(N, j-1) - w(N-1, j-1)\Big] \\ & - {}^{k^2 \lambda(N)}\!/\!{}_{Mh^3}\Big[w(N, j-1) - w(N-1, j-1)\Big]^3 \\ & + {}^{k^2}\!/\!{}_M\Big[d(j)\Big]\end{aligned} \tag{45}$$

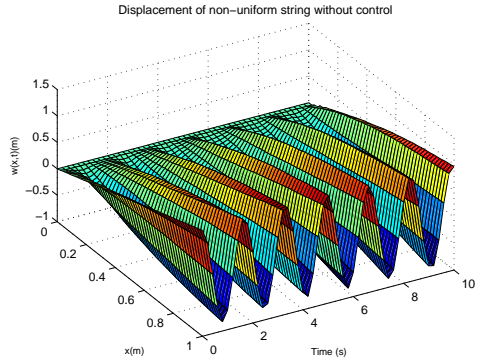

**Figure 7.** Non-uniform string - without control
$(T = 10000)$

**Remark.** *Non-uniform string without control are represented in Figure 7. Like above example, the string is still stable, although it undulates around 0 due to undesired external disturbance.*

**(ii) Exact model control**

$$\begin{aligned}u(j) =\ & T_0(L)w_x(L, j) - Mw_{xt}(L, j) - k_1 w_t(L, j) - \\ & k_2 w_x(L, j) - \mathrm{sgn}\left[w_t(L, j) + w_x(L, j)\right]\bar{d}\end{aligned} \tag{46}$$

The approximation for Eq. (37) after substituting Eq. (46)

$$\begin{aligned}w(N, j) =\ & 2w(N, j-1) - w(N, j-2) \\ & - {}^{Kk}\!/\!{}_M\Big[w(N, j-1) - w(N, j-2)\Big] \\ & + {}^{Kk^2}\!/\!{}_{Mh}\Big[w(N, j-1) - w(N-1, j-1)\Big] \\ & - {}^{k}\!/\!{}_h\Big[w(N, j-1) - w(N-1, j-1) \\ & \quad - w(N, j-2) + w(N-1, j-1)\Big] \\ & - {}^{2k^2}\!/\!{}_M\Big[sgn(A)\Big] + {}^{k^2}\!/\!{}_M\Big[d(j)\Big] \\ & - {}^{\lambda(N)k^2}\!/\!{}_{Mh^3}\Big[w(N, j-1) - w(N-1, j-1)\Big]^3 \\ & - {}^{10(L+1)k^2}\!/\!{}_{Mdx}\Big[w(N, j-1) \\ & \quad - w(N-1, j-1)\Big]\end{aligned} \tag{47}$$

#### 3.2.4 Step 4: Define general equation for each element

Finally, we obtain the governing equation from Eq. (34)

$$\begin{aligned}w(i, j) =\ & 2w(i, j-1) - w(i, j-2) \\ & + {}^{k^2 T(i,j)}\!/\!{}_{\rho h^2}\Big[w(i+1, j-1) - 2w(i, j-1) \\ & \quad + w(i-1, j-1)\Big] + {}^{k^2 T'(i,j)}\!/\!{}_{\rho h} \\ & \times \Big[w(i, j-1) - w(i-1, j-1)\Big] - {}^{k^2}\!/\!{}_\rho \Big[f(i, j)\Big] \\ & + {}^{k^2 \lambda'(i)}\!/\!{}_{\rho h^3}\Big[w(i, j-1) - w(i-1, j-1)\Big]^3 \\ & - {}^{3k^2 \lambda(i)}\!/\!{}_{\rho h^2}\Big[w(i, j-1) - w(i-1, j-1)\Big]^2 \\ & \times \Big[w(i+1, j-1) - 2w(i, j-1) + \\ & \quad w(i-1, j-1)\Big]\end{aligned} \tag{48}$$

#### 3.2.5 Step 5: Solve for derived or secondary quantities by recursion

Solution is obtained by using Matlab when we use all derived equations from former steps. System parameters and variables are defined in Table 2.

**Remark.** *Figure 8 and Figure 9 express the movement of the string with exact control while we choose suitable controller. The shape displacement goes down gradually before converging to stable position, 0, as we expect.*

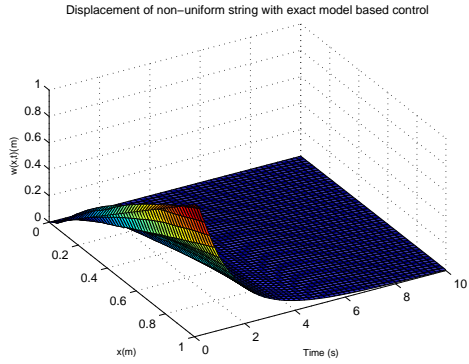

**Figure 8.** Non-uniform string - exact model control

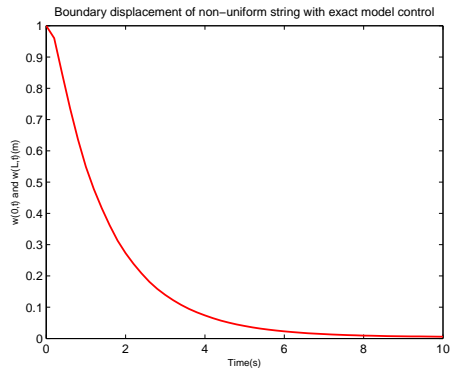

**Figure 9.** Non-uniform string - exact model control

#### 3.2.6  Step 6: Interpretation the result

Before coming out with the appropriate controller as well as mesh configurations like above, we tested the string with $T = 100$ and kept the left parameters.

**Remark.** *We can observe that instability of the string with $T = 100$ in Figure 10.*

## 4  SimuFlex

It is clear that many scientific & technology problems are governed by PDEs and hence this topic became an important and emerging field of research in the recent years. It has been observed that some of the problems have an analytical solution. PDEs can be classified as elliptic, parabolic or hyperbolic. The aim is to approximate solutions to differential equations that is to find a function (or some discrete approximation to this function) that satisfies a given relationship between its derivatives on some given region of space and time as well as some boundary conditions along the edges of this domain.

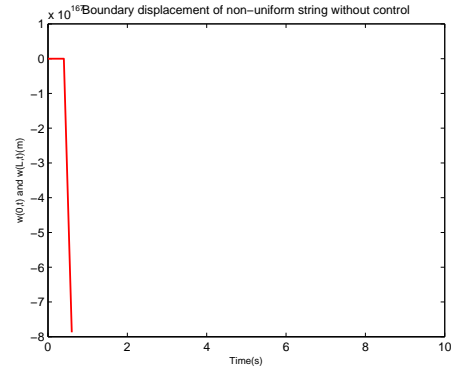

**Figure 10.** Non-uniform string - without control ($T = 100$)

Flexible systems including strings and beams are widely applied in many areas of modern mechanical fields. Elementary structural mechanics is primarily concerned with the behavior of line elements, such as rods or beams [19]. A beam is an element which carries load between the supports by virtue of its resistance to bending and shearing [20]. Defection computations may be used for computing the reactions for indeterminate beams and trusses as well as the stresses in the members of redundant trusses [21] [22].

The current practice of design for a given flexible system is troublesome. Engineers and researchers normally need to use a combination of tools for model building, controller design, and numerical simulation. Additional efforts are necessary to observe the simulation results to check whether the oscillation is damped out at the joints [23], [24]. Thus, a GUI was developed using MATLAB for the simulation of flexible beam in order to help users to easily use the simulation platform [25].

**Capabilities**

(i) Let the user input essential parameters such as approximation nodes, length, time, and controller gains.

(ii) Process data using C language and save results under MATLAB data files.

(iii) Read data and plot using MATLAB (3D for Displacement and 2D for Moving)

**Screenshots and Explanation**

The main features of this open-architecture platform in Figure 11 include:

(i) A list of systems which the user can choose from to simulate including Euler-Bernoulli beam, Timoshenko beam, exponential beam, string, and non-uniform string

(ii) A small window shows a figure, which usually describes a chosen system.

(iii) A group of control buttons: "Choose" button allows the user to take a general look at a system to simulate,

"Run" button runs the program. and "Close" button closes program.

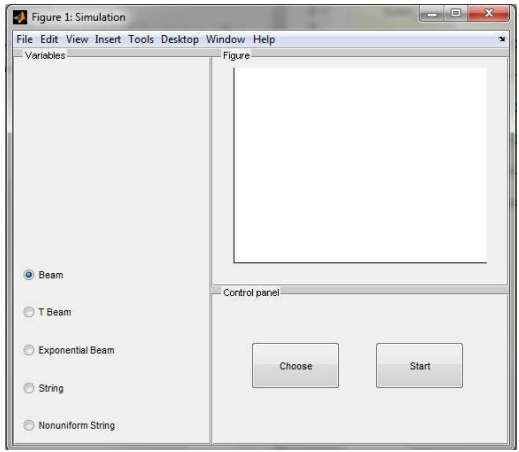

**Figure 11.** Starting page

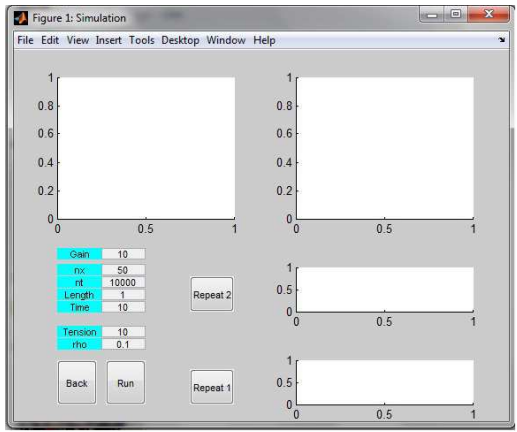

**Figure 12.** Simulation page - before running

The main features of the simulation platform is shown in Figure 12. Several workspaces can be opened and executed at the same time.

(i) A powerful system configuration is provided for easy configuration and viewing of various systems such as control gain, N, T, length, and time.

(ii) A flexible graphic plotting environment is embedded in the system for easily monitoring and observing the control performances. The user can view two working environments. Firstly, 3D figures help user view the overall performance of the system and check whether the designed system is stable or not. And then 2D moving figures help user to view the state of the flexible system over time.

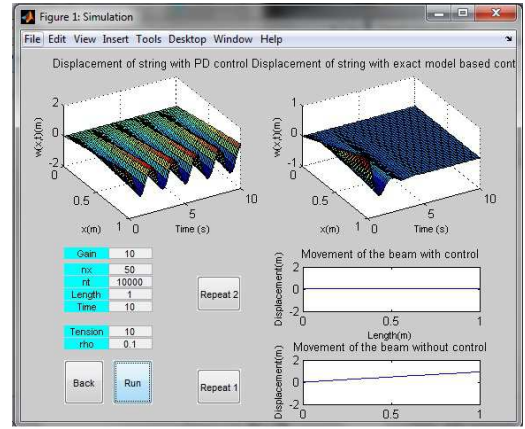

**Figure 13.** Simulation page - after running

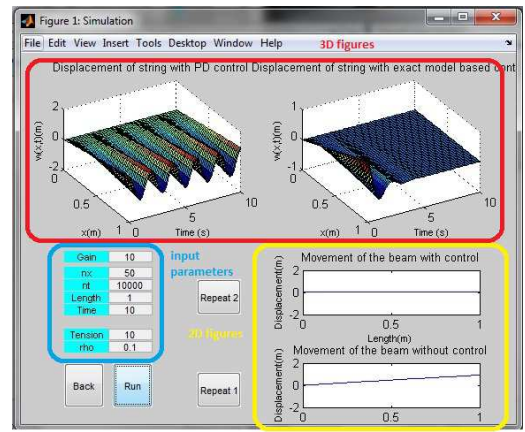

**Figure 14.** Simulation page - clear explanation

(iii) A group of control buttons: "Back" button returns to starting page, "Run" button runs the simulation, "Repeat 1" button replays the with control flexible system and "Repeat 2" button replays the without control flexible system.

## 5 Conclusions

In this article, a deep understanding of flexible systems and the way to model and simulate them under external disturbance has been studied. The main contribution of this article is building a general procedure for simulating flexible systems using finite difference method. In addition, contributions included formulating the general stabilizing condition for some specific systems with the purpose of verifying the effectiveness of the presented method. Next, two examples are provided to demonstrate the simulating procedure. Finally, a Graphical User Interface application, SimuFlex, for flexible systems is introduced.

Future works or applications of this proposed work can be, for example, applying the proposed procedure in a com-plete system, represented by partial differential equations with the purpose of simulation. Furthermore, researchers can find a general stabilizing condition for a system of partial difference equations.

# Biography

**Minh Hoang Vu**: He is an Erasmus Mundus scholar in Medical Imaging and Applications program. He has a B.Sc. in Control and Automation from National University of Singapore and a M.Sc. in Computer Control and Automation from Nanyang Technological University. His research interests are medical imaging, robotics, and machine vision. Email: hoangminh.vu@smartdatics.com

**Tajwar Abrar Aleef**: He received his Bachelor of Science in Electrical & Electronic Engineering from American International University-Bangladesh (2016). He is currently enrolled in Erasmus Mundus Joint Master Degree in Medical Imaging and Applications program between University of Burgundy (France), University of Cassino (Italy) and University of Girona (Spain). Email: tajwar_aleef@etu.u-bourgogne.fr

**Usama Pervaiz**: He is studying in Erasmus+ Joint Master Program in Medical Imaging and Applications; University of Burgundy (France), University of Cassino(Italy) and University of Girona (Spain). He has completed his bachelors in Electrical Engineering with major in signal and image processing from National University of Sciences and Technology, Pakistan. Email: 12beeupervaiz@seecs.edu.pk

**Yeman Brhane**: He received his BSc. in Electronics and Communication Engineering from University of Mekelle, Ethiopia. He is studying Erasmus Mundus Joint Master Degree in MedicAl Imaging and Applications: University of Buorgogne (France), University of Girona (Spain) and University of Cassino. His research interest includes 3D Image Registration, and Segmentation. Email: yemanbrhane1989@gmail.com

**Saed Khawaldeh**: He is an Erasmus Mundus scholar pursuing his masters in Medical Imaging and Applications program:University of Burgundy (France), University of Cassino (Italy) and University of Girona (Spain). He has B.Sc. in Computer and Electronics Engineering from Al-Balqa Applied University (Jordan). His research interests are Simultaneous EEG-fMRI Recordings, and Machine Learning. Email: khawaldeh.saed@gmail.com